*Universal Embedding spaces for G—manifolds*

Arthur G. Wasserman

*Abstract*  For any compact Lie group G and any n we construct a smooth G-manifold $U_n(G)$ such that any smooth n-dimensional G-manifold can be embedded in $U_n(G)$ with a trivial normal bundle. Furthermore, we show that such embeddings are unique up to equivariant isotopy It is shown that the (inverse limit) of the cohomology of such spaces gives rise to natural classes which are the analogue for G-manifolds of characteristic classes for ordinary manifolds. The cohomotopy groups of $U_n(G)$ are shown to be equal to equivariant bordism groups.

*Introduction*

If X and Y are smooth manifolds and $f : X \to Y$ is a smooth embedding then $f^* T(Y)$, the pullback of the tangent bundle of Y, can be written as $T(X) \oplus \nu(X, Y)$ where $\nu(X, Y)$ is the normal bundle of X in Y. If the dimension of Y is greater than twice the dimension of X then $\nu(X, Y)$ is in fact determined by the equation $T(X) \oplus \nu(X, Y) = f^* T(Y)$. [H]

Since $f^* T(Y)$ depends only on the homotopy class of f, $\nu(X, Y)$ is determined by the homotopy class of f and we will occasionally use the notation $\nu(f)$ rather than $\nu(X, Y)$ for the normal bundle when there is more than one embedding under consideration.

We now construct a manifold Y with the property that every n-dimensional manifold M admits an embedding into Y with a trivial normal bundle, that is, $\nu(M, Y)$ is isomorphic to $M \times \mathbb{R}^s$ for some s.

Here is a sketch of a proof of the existence of such a space that I believe was shown to me by Edgar. Brown Jr. many years ago.[B]

Let $s > 2n + 1$ and let $\mu_n \to G_n(\mathbb{R}^s)$ be the tautological n-dimensional vector bundle over $G_n(\mathbb{R}^s)$, the Grassmannian manifold of n-planes in $\mathbb{R}^s$, and let $\gamma \to G_n(\mathbb{R}^s)$ be a complement to the tangent bundle of $G_n(\mathbb{R}^s)$, that is, $\gamma \oplus T(G_n(\mathbb{R}^s))$ is a trivial bundle. [A]   Let Y = the total space of the bundle $\gamma \oplus \mu_n$ ; we claim that Y has the desired property.

Proof.  Let M be an n-manifold and let $\tau_M : M \to G_n(\mathbb{R}^s) \subset Y$ be the classifying map for the tangent bundle of M, that is, $\tau_M^*(\mu_n) = T(M)$. Since dimension $G_n(\mathbb{R}^s) = n(s-n)$ we can, for dimensional reasons, choose $\tau_M$ to be an embedding.

Next, $T(Y) |G_n(\mathbb{R}^s) = T(\gamma \oplus \mu_n) |G_n(\mathbb{R}^s) = T(G_n(\mathbb{R}^s)) \oplus T_F$ where $T_F$ is the tangent bundle along the fiber.  Since $T_F|G_n(\mathbb{R}^s) = \gamma \oplus \mu_n$ we have $T(Y) |G_n(\mathbb{R}^s) = T(G_n(\mathbb{R}^s)) \oplus \gamma \oplus \mu_n$ and thus
$T(M) \oplus \nu(M, Y) = \tau_M^* T(Y) = \tau_M^*(\mu_n) \oplus \tau_M^*(T(G_n(\mathbb{R}^s)) \oplus \gamma)$. By definition $\tau_M^*(\mu_n) = T(M)$ and $T(G_n(\mathbb{R}^s)) \oplus \gamma$ is a trivial bundle so $\tau_M^*(T(G_n(\mathbb{R}^s)) \oplus \gamma)$ is a trivial bundle.
Hence $T(M) \oplus \nu(M, Y) = T(M) \oplus$ a trivial bundle and so $\nu(M, Y)$ is a trivial bundle since we are in the stable range. [H, page 112]  Thus M embeds in Y with trivial normal bundle.

Note that Y is <u>not</u> unique; for example, any n-manifold can also be embedded in $Y \times \mathbb{R}$ or $Y \times S^{n+1}$.

We can easily construct a compact embedding space namely $D(\gamma \oplus \mu_n) \subset Y$, the disk bundle of $\gamma \oplus \mu_n$. In fact, since the codimension of the zero section of $D(\gamma \oplus \mu_n)$ is greater than n we can assume by transversality that the image of any embedding is disjoint from the zero section and thus retracts to $S(\gamma \oplus \mu_n)$. In other words we can construct an embedding space that is a closed manifold.

In this paper the above result is generalized to show the existence of such embedding spaces for G-manifolds, that is, for manifolds with a group action.

In section 1 we prove the existence of embedding spaces for several special cases of G-manifolds, namely tubular neighborhoods of a stratum in a G-space, and then, in section 2, combine those cases into the final result.

In section 3 we discuss applications and list some open questions.

§1 Embedding spaces for special G-manifolds

Preliminary Definition: A manifold $Y_n$ is an n-universal embedding space if every n-dimensional manifold M admits an embedding into $Y_n$ with a trivial normal bundle, that is, $\nu(M, Y_n)$ is isomorphic to $M \times \mathbb{R}^s$ for some s and furthermore, any two such embeddings are isotopic via an isotopy having a trivial normal bundle, that is, $\nu(M \times I, Y_n \times I)$ is isomorphic to $M \times I \times \mathbb{R}^s$ for some s.

The defnition will be slightly modified in section 3.

Let $s > 2n + 1$ and let $\mu_{n+1} \to G_{n+1}(\mathbb{R}^s)$ be the tautological (n+1)-dimensional vector bundle over $G_{n+1}(\mathbb{R}^s)$, the Grassmannian manifold of (n+1)-planes in $\mathbb{R}^s$, and let $\gamma \to G_{n+1}(\mathbb{R}^s)$ be a complement to the tangent bundle of $G_{n+1}(\mathbb{R}^s)$, that is, $\gamma \oplus T(G_{n+1}(\mathbb{R}^s))$ is a trivial bundle. [A]

We note that $T(Y_n)|G_{n+1}(\mathbb{R}^s) = T(G_{n+1}(\mathbb{R}^s)) \oplus T_F = T(G_{n+1}(\mathbb{R}^s)) \oplus \gamma \oplus \mu_{n+1} = \mu_{n+1} \oplus$ a trivial bundle.

Let $Y_n$ = the total space of $\gamma \oplus \mu_{n+1}$ as above or $D(\gamma \oplus \mu_{n+1})$ or $S(\gamma \oplus \mu_{n+1})$; we claim that $Y_n$ has the desired property.

Moreover, we can slightly extend the observations of the introduction to the following proposition.

Proposition 1. If $f : M^n \to A$ is smooth then there is an embedding $h : M^n \to Y_n$, unique up to isotopy, such that $\nu(f \times h)$ is trivial where $f \times h : M^n \to A \times Y_n$.

Proof: For any choice of h we have $T(M) \oplus \nu(M, A \times Y_n) = (f \times h)^* T(A \times Y_n) = f^* T(A) \oplus h^* T(Y_n) = f^* T(A) \oplus h^*(\mu_{n+1}) \oplus$ a trivial bundle.

Let $\xi \to M$ be a vector bundle over M such that $\xi \oplus f^*(T(A))$ is trivial. Since dimenension M = n we can write $\xi \oplus T(M)$ uniquely as $\zeta \oplus$ a trivial bundle where dimension $\zeta = n+1$.[A].

Choose h such that $h^*(\mu_n) = \zeta$. ( Such an h is unique up to homotopy and hence isotopy for dimensional reasons.) Then $f^* T(A) \oplus h^* T(Y_n) = f^* T(A) \oplus \xi \oplus T(M) \oplus$ a trivial bundle = $T(M) \oplus$ a trivial bundle and hence $T(M) \oplus \nu(M, A \times Y_n) = T(M) \oplus$ a trivial bundle so $\nu(M, A \times Y_n) = \nu(M, f \times h)$ is trivial since we are in the stable range.

Note 1: We need $\mu_{n+1} \to G_{n+1}(\mathbb{R}^s)$ rather than $\mu_n \to G_n(\mathbb{R}^s)$ to get an h that is unique up to isotopy.

Note 2: the proof of proposition 1 uses one of the specific n-universal embedding spaces $Y_n$ above and we will have occasion to use a specific n-universal embedding space $Y_n$ in later propositions. (It is not clear if the proposition is true for an arbitrary n-universal embedding space. )

Notation: Let $X_n$ denote the n-universal embedding space for unoriented manifolds given above by $S(\gamma \oplus \mu_{n+1})$.

One could similarly construct n-universal embedding space for other types of manifolds, for example, for orientable manifolds using the Grassmannian manifold of oriented (n+1)-planes in $\mathbb{R}^s$; proposition 1 would still be true for A an orientable manifold. However, in this paper we will only consider the unoriented case.

Proposition 2 If $\pi: E \to B$ and $\pi': E' \to B'$ are fiber bundles in the sense of Steenrod [St] with fiber F and $\Phi: E \to E'$ is a bundle map then $\nu(E, E') = \pi^* \nu(B, B')$.

Proof. By definition $T(E) \oplus \nu(E, E') = \Phi^* T(E')$. We have $T(E) = \pi^* T(B) \oplus T_F$ and $T(E') = \pi'^* T(B') \oplus T'_F$ where $T_F$ is the tangent bundle along the fiber of $\pi: E \to B$ and $T'_F$ is the tangent bundle along the fiber of $\pi': E' \to B'$.

Let $\phi: B \to B'$ be the map induced on the base spaces by $\Phi$; then $\pi' \Phi = \phi \pi$ and hence, $(\pi' \Phi)^* T(B') = (\phi \pi)^* T(B') = \pi^* \phi^* T(B')$.

Thus, $\pi^* T(B) \oplus T_F \oplus \nu(E, E') = \Phi^* T(E') = \Phi^*(\pi'^* T(B') \oplus T'_F) = \Phi^*(\pi'^* T(B')) \oplus \Phi^* T'_F$. Since $\Phi$ is a bundle map $\Phi^* T'_F = T_F$ thus $\pi^* T(B) \oplus \nu(E, E') = \Phi^*(\pi'^* T(B')) = \pi^* \phi^* T(B')$. By definition $T(B) \oplus \nu(B, B') = \phi^* T(B')$ so $\pi^* T(B) \oplus \pi^* \nu(B, B') = \Phi^*(\pi'^* T(B'))$ and the result follows.

We now consider the existence of n-universal embedding spaces for G-spaces, that is, spaces with an action of a compact Lie group G.

Notation: If X is a G-space then $G_x = \{g \in G \mid gx = x\}$ is the isotropy subgroup at x. If $H \subset G$ is a subgroup then $X^H = \{x \in X \mid H \subset G_x\} = \{x \in X \mid hx = x \text{ for all } h \in H\}$ is the fixed point set of H.

We start with a few special cases that illustrate the method before dealing with the details of the general case

First we construct an n-universal embedding space for free actions of a compact Lie group, that is, G-manifolds M such that the isotropy subgroup at x, $G_x$ = e for every x ∈ M.

So let G be a compact Lie group and let EG(n) be an n-universal bundle for G in the sense of Steenrod , that is, if G acts freely on an G-space P with quotient space K, an n-dimensional simplicial complex, Q is a subbundle of P with quotient L , a subcomplex of K, and h : Q → EG(n) is an equivariant map then there is an equivariant map f from P to EG(n) such that f|Q = h. The existence of such a space was shown in Steenrod [St, §19.6]. The Steenrod construction EG((n) is, in fact, a closed smooth G-manifold so that if h is a smooth equivarinat map then we may choose the extension f to also be be a smooth extension [W Cor 1.12 ].

Note: The use of universal in this statement differs from its use in the phrase "n-universal embedding space".

Lemma 3 Let $M^n$ be a free G-manifold and let f: $M^n$ → C be a smooth equivariant map of G-manifolds . Let m = dimension M - dimension G; then there is an equivariant embedding h: $M^n$ → EG(m + 1) x $X_m$ , unique up to isotopy, such that ν(f × h ) is trivial where f × h: $M^n$ → C × EG(m + 1) × $X_m$.

Proof: Since EG(m+1) is an (m + 1)-universal bundle for G there is a smooth equivariant map
k: $M^n$ → EG(m + 1 ) and hence a smooth equivariant map f × k : $M^n$ → A × EG(m + 1) which induces a smooth map of quotient spaces (f × k)' : M/G → (A × EG(m + 1))/G. (The quotient spaces are manifolds because G acts freely). We now apply proposition 1 to the map (f × k)' noting that
dimension M/G = dimension M - dimension G = m [P 1.7.31] to get a map j : M/G → $X_m$ such that ν( (f × k)' × j) is trivial where (f × k') × j: $M^n$/G → (X× EG(m + 1))/G × $X_m$ . Let h : M → EG(m + 1) × $X_m$ be given by k × j. Finally, applying proposition 2 to the bundle maps
π : M→ M/G and π' :A × EG(m + 1) × $X_m$→ (A × EG(m + 1) )/G x $X_m$ we get that ν( f × h ) is trivial as we wished to show. Uniqueness follows as in Lemma 1.

Remark 1 Taking C = point in Lemma 3 shows that EG(m + 1) x $X_m$ is an n-universal embedding space for n-dimensional free G manifolds.

We next fix a conjugacy class of subgroups (H) and consider G-manifolds M with $G_x$ conjugate to H ⊂ G for every x ∈ M.

Some notation: If H ⊂ G and X is an H space then Y = G $\times_H$ X = G × X/{(gh, x) = (g, hx)} is a G-space with G action given by g'(g, x) = (g'g, x). Note that Y/G = X/H. Moreover, if f: X → Z is an H-equivariant map then f extends to a G-equivariant map F: G $\times_H$ X → G $\times_H$ Z. [P]

Let N = N(H) be the normalizer of H in G, N' = N(H)/H, $M_H$ = {x ∈ M | $G_x$ = H }; then $M_H$ is an N'-space. If $M_{(H)}$ = {x ∈ M | $G_x$ is conjugate to H} then $M_{(H)}$ = G $\times_N$ $M_H$ = $GM_H$. That is, M fibers over

G/N with fiber $M_H$; thus dimension $M$ = dimension $M_H$ + dimension G/N.

**Lemma 4** Let $f : M^n \to C$ be a smooth equivariant map of G-manifolds and suppose that $M = M_{(H)}$. Let $s \geq$ dimension M - dimension G/H ; then there is an equivariant embedding
h: $M^n \to G \times_N (C^H \times EN'(s + 1)) \times X_s$ , unique up to equivariant isotopy, such that
$\nu(h)$ is trivial and the composition h: $M^n \to G \times_N (C^H \times EN'(s + 1)) \times X_s \to G \times_N C^H \to C$ is f.

Proof: Since N' acts freely on $M_H$ and $f(M_H) \subset C^H$ we apply Lemma 3 (with N' playing the role of G) to $f|M_H : M_H \to C^H$ to get a map k: $M_H \to EN'(s + 1) \times X_s$ such that $\nu(f|M_H \times k)$ a trivial bundle. Extending the map $f|M_H \times k$ to $M = G \times_N M_H$ we get the desired map h via
$M = G \times_N M_H \to G \times_N (C^H \times EN'(s + 1) \times X_s)$. To see that $\nu(h)$ is a trivial bundle we pass to the quotient manifolds getting $(G \times_N M_H)/G = M_H/N \to (G \times_N C^H \times EN'(s + 1) \times X_s)/G =$
$(C^H \times EN'(s + 1))/N \times X_s$ which is an embedding with trivial normal bundle by construction. Taking the quotient by G is a fibration with fiber G/H ([P]) and hence, by Proposition 2, $\nu(h)$ is a trivial bundle. The other claims are clear.

Remark 2 Taking C = point in Lemma 4 shows that $G \times_N EN'(s + 1) \times X_s$ is an n-universal embedding space for n-dimensional G manifolds having only one conjugacy class of isotropy subgroups.

Some notation: Let V be a representation of the group G. We denote by $G_n(V)$, the Grassmannian manifold of n-planes in V with the obvious G action. Let $\mu_n(V) \to G_n(V)$ denote the tautological n-dimensional G-vector bundle over $G_n(V)$.

If V is a representation of G then t V denotes $V \oplus V \oplus \ldots V$, t summands.

We recall a definition and theorem from [W1].

Definition. Let $\pi : E \to M$ be a G-vector bundle of fibre dimension $k < \infty$ over the G-manifold M. The bundle is said to be subordinate to the representation V of G if, for each $m \in M$, the representation of $G_m$ on $E_m = \pi^{-1}(m)$ is equivalent to a subrepresentation of t $V|G_m$ for some t.

Remark 3 It follows immediately from the definition that if $X = A \cup B$ and $E \to X$ is a G-vector bundle with E|A and E|B both subordinate to a representation V of G then E is subordinate to V.

Theorem 5 [W1] The set of equivalence classes of k-dimensional G-vector bundles over $M^n$ subordinate to V is isomorphic to the set of equivariant homotopy classes of maps of M into $G_k(tV)$ if $t > n + k + 1$.

Notation The case $V = \mathbb{R}$ is special; we write $G_k(\mathbb{R}^t)$ rather than the consistent but strange looking $G_k(t\mathbb{R})$.

Remark 4 Note that an equivariant map **from** M to $G_k(\mathbb{R}^t)$ must factor through M/G since G acts trivially on $G_k(\mathbb{R}^t)$. Thus a bundle over M subordinate to $\mathbb{R}$ is the pullback of a bundle over M/G.

We next construct n-universal embedding spaces for n-dimensional G-manifolds E which are total spaces of k-dimensional G-vector bundles subordinate to a representation V of G and with base space a G-manifold $M^{n-k}$.

Lemma 6  Let $M^{n-k}$ be a G-manifold with $M = M^G$ and let $\pi : E \to M$ be a k-dimensional G-vector bundle subordinate to a representation V of G with $V^G = 0$. If $t > n+1$ then there is an equivariant embedding of E into $\mu_k(tV) |G_k(tV)^G \times X_{n-k}$ with a trivial normal bundle; the embedding is unique up to equivariant isotopy.

Proof: Let $j: M \to G_k(tV)$ be the classifying map for the G-vector bundle E and let $J : E \to \mu_k$ be the covering bundle map ; the existence of such a map is provided by theorem 5 since $t > (n-k) + k + 1 = n + 1$.   Since $M = M^G$, $j(M) \subset G_k(tV)^G$ so let j' denote the map $j|M : M \to G_k(tV)^G$ and let $J' : E \to \mu_k |G_k(tV)^G$ be the covering bundle map. By Proposition 1 there is a map $h : M \to X_{n-k}$ such that $\nu(j' \times h)$ is trivial. By Proposition 2  $\nu( J' \times h)$ is the pull back of the trivial bundle $\nu(j' \times h)$ and hence is also trivial. That proves the existence of an equivariant embedding $J' \times h$ with trivial normal bundle.

To show uniqueness of the map $J' \times h$ we let  $a \times b: E \to \mu_k |G_k(tV)^G \times X_{n-k}$ be another such equivariant embedding with trivial normal bundle. We note that for each $x \in M$, $E_x^G = 0$ since $V^G = 0$ so $M = E^G$. By equivariance $a \times b$ maps fixed point sets to fixed point sets so $a \times b$ maps M to $G_k(tV)^G \times X_{n-k}$. We have that $T(E) |M = T(M) \oplus E$, $T( \mu_k |G_k(tV)) = T(G_k(tV)^G ) \oplus \mu_k |G_k(tV)^G$  and, for $y \in E$, $p(y) = x \in M$, $da_x(y) \in \mu_k |G_k(tV)^G$ , that is, $da_x$ respects the direct sum decomposition  because $da_x$ is an equivariant linear map and $E_x^G = 0$ so Hom $(E_x, T(G_k(tV)a_{(x)})^G) = 0$.

We next linearize the map $a \times b$ by the isotopy $f_t(y) = ( \frac{1}{t} a (ty) , b(ty))$ for $t > 0$ and $f_0(y) = (da_{p(y)}(y) , b \circ p(y) )$ where p is the projection $E \to M$ and $da_{p(y)}(y)$ denotes as usual the differential of the map a at the point p(y) in M.  As noted above, $f_o$ is an equivariant bundle map and hence by Theorem 5 $f_o|M$ is equivariantly homotopic to $j \times b$ .

Since f is equivariantly isotopic to the bundle map $f_0$ the projection of $f_0$ to $G_k(tV)^G$ is equivariantly isotopic to the classifing map j' and an application of lemma 4 in the special case $H = G$ yields the uniqueness result.

We now give a simple example to show that the condition $V^G = 0$ in Lemma 6 is necessary.

Example: Let $G = Z_2$ and let $V = \mathbb{R} \oplus \mathbb{R}'$ where G acts trivially on $\mathbb{R}$ and nontrivially on $\mathbb{R}'$ , let L denote the nontrivial line bundle over $S^1$ with trivial G action, let L' denote the nontrivial line bundle over $S^1$ with G action given by $g x = -x$ and let $S^1 \times \mathbb{R}$ denote the trivial line bundle over $S^1$ with trivial G action. Let $E = L \oplus L' \oplus \mathbb{R}$. We may regard E as a G-vector bundle over L with fiber $\mathbb{R} \oplus \mathbb{R}'$ and also as a G-vector bundle over $S^1 \times \mathbb{R}$ again with fiber $\mathbb{R} \oplus \mathbb{R}'$. As in lemma 6 there is a classifying map for the bundle $E \to L$ into $G_2(tV)^G$ specifically into the component diffeomorphic to $G_1(\mathbb{R}^t) \times G_1(t\mathbb{R}')$
and the map $(a,b) : L \to G_1(\mathbb{R}^t) \times G_1(t\mathbb{R}')$ classifies $S^1 \times \mathbb{R}$ via a and L' via b so a is homotopically trivial and b is not. However, if we regard E as a G-vector bundle over $S^1 \times \mathbb{R}$ the classifying map (a', b') : $S^1 \times \mathbb{R} \to G_1(\mathbb{R}^t) \times G_1(t\mathbb{R}')$ classifies L via a' and L' via b so a' is homotopically nontrivial and b' is not. The resulting equivariant embeddings E $\to G_2(tV)^G \times X_2$ are clearly not homotopic thus not isotopic and the conclusion of lemma 6 is false without the assumption that $V^G = 0$.

The n-universal embedding spaces for n-dimensional G-manifolds E which are total spaces of k-dimensional G-vector bundles subordinate to a representation V of G and with base space a G-manifold $M^{n-k}$ constructed above are not connected in general. The components of the spaces $G_k(tV)^G$ can be indexed by the k-dimensional representations of G subordinate to V [c.f. [DM], [DW]]. If W is such a representation—which just means every irreducible subrepresentation of W is also a subrepresentation of V—then we call the component of $G_k(tV)^G$ associated to W, $G_W(tV)^G$. It is easy to see that $\mu_k | G_W(tV)^G \times X_{n-k}$, a component of $\mu_k | G_k(tV)^G \times X_{n-k}$, is an n-universal embedding space for G-manifolds E which are total spaces of G-vector bundles with $E_x = W$ as a representation of G for all x in M. However, the manifolds $G_W(tV)^G$ vary in dimension so we must "pad" them out by multiplying by some factor $\mathbb{R}^s$ where s depends on W to have a manifold in the usual sense.

That is, we fix t = n + 2 **and** note that dimension $\mu_k | G_W(tV)^G \leq k + k(tk - k) \leq n^2 + n$ so we may chose s = dimension $X_n + n^2 + n -$ dimension $\mu_k | G_W(tV)^G \times X_{n-k}$.

Thus each component will be of the form $\mu_k | G_W(tV)^G \oplus \mathbb{R}^s \times X_{n-k}$. We will further "pad" each component so that the dimensions are the same for all k as well as for all W.

Definition: Let $U_n(G,W) = \mu_k | G_W((n+2)V)^G \oplus \mathbb{R}^s \times X_{n-k}$
where V is some fixed representation of G such that W is subordinate to V and dimension $W \leq n$.

We next construct n-universal embedding spaces for n-dimensional G-manifolds E which are total spaces of k-dimensional G-vector bundles with base space a G-manifold $M^{n-k} = E^G$.

Let $O_n(G,k)$ = the disjoint union over the distinct k-dimensional representations W of G of the $U_n(G,W)$ and let $O_n(G)$ = the disjoint union of the $O_n(G,k)$ for k = 1, 2, ...n.

It is not true that any such bundle is subordinate to a representation of G. For example, let $\mathbb{C}_1, \mathbb{C}_2, ..$ denote the irreducible 2 dimensional representations of G = $S^1$ and let E = the disjoint union of the $\mathbb{C}_j$'s and let $\pi$: E $\to \mathbb{Z}_+$ be given by $\pi(\mathbb{C}_j) = \{j\}$. E is obviously not subordinate to any finite dimensional representation of G. However, if M1 is a connected component of M then E|M1 is trivially subordinate to the representation of G on $E_x$ where x is any point in M1. Thus we have:

Lemma 7 Let $M^{n-k}$ be a G-manifold and let $\pi : E \to M$ be a k-dimensional G-vector bundle with $M^{n-k} = E^G$. Then there is an equivariant embedding of E into $O_n(G)$ with a trivial normal bundle; the embedding is unique up to equivariant isotopy.

More generally we have:

Proposition 8  Let $A^n$ be a G-manifold. Then there is an invariant neighborhood of $A^G$ that embeds into $O_n(G)$ with a trivial normal bundle; the embedding is unique up to equivariant isotopy.

Proof $A^G$ is a disjoint union of closed subnmanifolds $A_0, A_1, A_2, ...A_n$ of codimensions 0, 1, 2 ..n respectively and each $A_j$ has a tubular neighborhood that embedds in $O_n(G,k)$ hence $A^G$ has an invariant neighborhood that embedds into $O_n(G)$ with a trivial normal bundle; the embedding is unique up to equivariant isotopy.

We now fix a normal subgroup $H \subset G$ and a representation V of G and consider manifolds E which are the total spaces of k-dimensional G-vector bundles subordinate to V over G-manifolds $M = M_H$ and construct n-universal embedding spaces for such G-manifolds.

Lemma 9  Let $M^{n-k}$ be a G-manifold with $M = M_H$ and let $\pi : E \to M$ be a k-dimensional G-vector bundle subordinate to a representation V of G with $V^G = 0$. Then if $t > n+1$ there is an equivariant embedding of E into $\mu_k |G_k(tV)^H \times E(G/H)(s + 1)) \times X_s$ with a trivial normal bundle where

$s = n - k - \dim G/H$.

Proof: Let $f : M \to G_k(tV)$ be the classifying map for the G-vector bundle E. The existence of such a map is provided by theorem 5 since $t > (n - k) + k + 1 = n + 1$. Since $M = M_H$, image $(f) \subset G_k(tV)^H$ so let $f' : M \to G_k(tV)^H$. Then $f'$ is covered by an equivariant bundle map $F: E \to \mu_k|G_k(tV)^H$. Now G/H acts freely on M so we may apply lemma 3 to $f'$ to get an equivariant embedding
$h : M^{n-k} \to E(G/H)(s+1)) \times X_s$, unique up to equivariant isotopy, such that $\nu(f' \times h)$ is trivial. Thus, by proposition 2, $F \times (h \circ \pi) : E \to \mu_k|G_k(tV)^H \times E(G/H)(s+1)) \times X_s$ is an equivariant embedding of E in $\mu_k|G_k(tV)^H \times E(G/H)(s+1)) \times X_s$ with a trivial normal bundle.

We now fix a conjugacy class of subgroups (H) of G and a representation V of G and consider manifolds E which are the total spaces of k-dimensional G-vector bundles subordinate to V over G-manifolds $M = M_{(H)}$ and construct n-universal embedding spaces for such G-manifolds. As usual N = normalizer of H in G and N' = N/H.

Proposition 10 Let $M^{n-k}$ be a G-manifold with $M = M_{(H)}$ and let $\pi : E \to M$ be a k-dimensional G-vector bundle subordinate to a representation V of G with $V^G = 0$. Then if $t > n + 1$ there is an equivariant embedding f of E into $G \times_N (\mu_k|G_k(tV)^H \times EN'(s+1)) \times X_s$ with a trivial normal bundle unique up to isotopy where $s = n - k - \dim G/H$.

Proof: Lemma 9 applies to the k-dimensional N(H) vector bundle $E|M_H \to M_H$ so let
$F \times (h \circ \pi) : E|M_H \to \mu_k|G_k(tV)^H \times EN'(s+1)) \times X_s$ be the N(H) equivariant embedding of $E|M_H$ in $\mu_k|G_k(tV)^H \times E(N/H)(s+1)) \times X_s$ with a trivial normal bundle that was constructed in the proof of lemma 9. We extend $F \times (h \circ \pi)$ by G-equivariance to
$G \times_N E|M_H \to G \times_N (\mu_k|G_k(tV)^H \times EN'(s+1)) \times X_s$.
We note that $M = M_{(H)} = G \times_N M_H$ and $E = G \times_N E|M_H$ thus we have the desired map
$E \to G \times_N \mu_k |G_k(tV)H \times E(N')(s+1)) \times X_s$.

To verify that this map is indeed an embedding with trivial normal bundle it is sufficent by proposition 2 to verify that
$G \times_N M_H \to G \times_N G_k(tV)^H \times E(N')(s+1)) \times X_s$ is an embedding with trivial normal bundle.

Passing to the quotient space we have $M_H/N' \to G_k(tV)^H \times_{N'} E(N')(s+1) \times X_s$ which is an embedding by the construction in lemma 9.

Consider now pairs (H, W) where H is a subgroup of G and W is a representation of H. We say that (H, W) ~ (H', W') if $G \times_H W$ is equivariantly diffeomorphic to $G \times_{H'} W'$. We denote the equivalence class of (H,W) by [H,W].

We say that a G-vector bundle $E \to M$ is of type $[H,W]$ if for each $x \in M$, $(G_x, E_x) \in [H, W]$.

**Proposition 11** If $W^H = 0$ there is a G-manifold $O_n(G, [H,W])$ such that any bundle of type $[H,W]$ embeds in $O_n(G, [H,W])$ with trivial normal bundle and such embeddings are unique up to equivariant isotopy.

Proof Let $E \to M$ be any bundle of type $[H,W]$. We cannot quote proposition 10 directly because we do not have a respresentation V of G such that E is subordinate to V. However, if H has finite index in G, for example if G is finite or H contains the identity component of G, we can let $V = i_H^* W$, the induced representation. By Frobenius reciprocity, $V|H$ contains every irreducible representation of H that W contains and $V^G = W^H = 0$ so $E \to M$ satisfies the conditions of proposition 10.
If H does not have finite index in G then V is still defined, see [Bo], but it is not finite dimensional. The induced representation does still obey Frobenius reciprocity in this case [S] and it is decomposable into finite dimensional representations of G so if W1 is an irreducible summand of W and V1 is an irreducible summand of $i_H^* W_1$ then $V_1|H$ will contain $W_1$ so putting $V = \oplus V_j$ where j runs over the irreducible components of W yields the desired representation of G and we can apply proposition 10 to this case also.

Proposition 11 includes as special cases proposition 1 and 10 and lemmas 3,4,6, 7, 8 and 9. In the next section we put togther the universal embedding spaces constructed for various strata to get universal embedding spaces for arbitrary G-manifolds.

## §2 Embedding spaces for G-manifolds

In this section we construct n-universal embedding spaces for general G-manifolds.

Let G be a compact Lie group and let $\mathcal{O}$ be a set of pairs (H,V) where H is a closed subgroup of G and V is a representation of H.

We say that $\mathcal{O}$ is admissible if

i) for any $(H,V) \in \mathcal{O}$ and any $x \in G \times_H V$ the pair $(G_x, S_x) \in \mathcal{O}$ where $G_x$ is the isotropy subgroup of x and $S_x$ is the slice representation of $G_x$ at x ( [K] ) and

ii) $(H,V) \in \mathcal{O}$ if and only if $(H, V \oplus \mathbb{R}) \in \mathcal{O}$.

Let $\mathcal{C}(G, \mathcal{O})$ be the category of G-manifolds M with $(G_x, S_x) \in \mathcal{O}$ for every $x \in M$ and morphisms being equivariant embeddings with trivial normal bundles.

Categories of manifolds defined by orbit types have been used in [J] and [LW].

Condition i) is just a consistency condition guaranteeing that $G \times_H V \in \mathcal{C}(G, \mathcal{O})$
if and only if $(H,V) \in \mathcal{O}$.
Condition ii) guarantees that $M \times I \in$

$\mathfrak{c}(G, \mathcal{O})$ and $\partial M \in \mathfrak{c}(G, \mathcal{O})$ if $M \in \mathfrak{c}(G, \mathcal{O})$.

If $M \in \mathfrak{c}(G, \mathcal{O})$ then any submanifold of M with trivial normal bundle $\in \mathfrak{c}(G, \mathcal{O})$.

As defined previously, we say that $(H, W) \sim (H', W')$ if $G \times_H W$ is equivariantly diffeomorphic to $G \times_{H'} W'$. We denote the equivalence class of (H,W) by [H,W]. Let $\mathcal{O}' = \{[H,W] \mid (H,W) \in \mathcal{O}\}$.

Define $|\mathcal{O}'|$ to be the number of $[H,V] \in \mathcal{O}'$ such that $V^H = 0$.

Our first objective is to prove the following:

Theorem 14 If $\mathcal{O}$ is admissible and $|\mathcal{O}'|$ is finite then there exists an n-universal embedding space $U_n(G, \mathcal{O}) \in \mathfrak{c}(G, \mathcal{O})$ that is, any $M^n \in \mathfrak{c}(G, \mathcal{O})$ admits an equivariant embedding in $U_n(G, \mathcal{O})$ with a trivial normal bundle and furthermore any two such embeddings are isotopic.

Before we prove the above theorem we note some facts about transversality in $\mathfrak{c}(G, \mathcal{O})$. Transversality for equivariant maps can be a delicate issue but the following facts are routine applications of Sard's theorem and partitions of unity.

a) If $f: M \to \mathbb{R}^s$ is a smooth equivariant map then f can be approximated arbitrarily closely by a smooth equivariant map f' transverse to $\{0\}$ even in the Whitney topology.

b) If $f: M \to W$ is a smooth equivariant map and $B \subset W$ is a closed submanifold with $\nu(B, W)$ trivial then f can be approximated arbitrarily closely by a smooth equivariant map f' transverse to B even in the Whitney topology.

c) If $E \to X$ is a G-vector bundle subordinate to $\mathbb{R}$ then any equivariant section of E can approximated arbitrarily closely by a smooth equivariant section transverse to the zero section. even in the Whitney topology.

d) If $f: M \to W$ is a smooth equivariant map and $B \subset W$ is a closed submanifold with $\nu(B, W)$ subordinate to $\mathbb{R}$ then f can be approximated arbitrarily closely even in the Whitney topology by a smooth equivariant map f' transverse to B.

Lemma 12 Let f: $M^n \to A$ be an equivariant embedding of G-manifolds with normal bundle subordinate to $\mathbb{R}$; then there is a smooth equivariant map $h: M \to X_n$, unique up to homotopy, such that $f * h: M \to A \times X_n$ has an equivariantly trivial normal bundle.

Proof We have $f^*(T(A)) = T(M) + \nu(M, A)$ where $\nu(M, A)$ is subordinate to $\mathbb{R}$. Let $\xi^{n+1}$ be the unique (n+1)-dimensional bundle such that $\nu(M, A) \oplus \xi^{n+1}$ is a trivial bundle. Such a bundle exists since $\nu(M, A)$ is classified **by** a map **to** $G_{n+1}(\mathbb{R}^{2n+1})$ **and by** remark d) about transversality above. (It is a trival extension of [A] **and** [H p.125].)

Let h' : M → $G_{n+1}(\mathbb{R}^{2n+1})$ be the classifying map **for** $\xi^{n+1}$ **and** let h : M → $X_n$ be chosen so that h' = π ∘ h where π : $X_n$ → $G_{n+1}(\mathbb{R}^{2n+1})$.

Since $T(X_n) = \mu_{n+1}(\mathbb{R}^{2n+1})$ ⊕ trivial we have h*(T($X_n$)) ⊕ ν(M,A) is trivial.

To complete the proof we first note that if i: A → B and j : B → C are embeddings then ν(A, C) = ν(A, B) ⊕ i*ν(B, C), that is, ν(j ∘ i) = ν(i) ⊕ i*ν(j).

Let U be an invariant tubular neighborhood of M in A with projection p : U → M and let g : U → U × $X_n$ be given by x → (x, h(p(x)). Then the composition g ∘ f = h. Since ν(g ∘ f) = ν(f) ⊕ f*ν(g) = ν(f) ⊕ f*ν(g ∘ f) = ν(f) ⊕ h*ν(g) = ν(M, A) ⊕ h*(T($X_n$)) = trivial bundle.

The extended argument given above was necessary because in general there is no stable range for G-vector bundles.

The following proposition is for illustrative purposes only; it demonstrates the ideas of the proof of the main theorem without some of the obscuring complications of the general case.

Let G be a compact Lie group and let V be a k-dimensional representation of G such that $G_x$ = e for all x ∈ V-{0}. For the next proposition we consider G-manifolds $M^n$ such that for all x ∈ M , $G_x$ = e or $G_x$ = G and furthermore, ν($M^G$, M) is a k-dimensional G-vector bundle subordinate to V. That is, we let $\mathcal{O}$ be the smallest admissible set of orbit types containing { (e, $\mathbb{R}$), (G, $W_i$)} where the $W_i$ are all the k-dimensional representations of G whose irreducible factors are also irreducible factors of V.

Proposition 13 There is an n-universal embedding space $U_n$ (C(G, $\mathcal{O}$)) for G manifolds in C(G, $\mathcal{O}$).

Proof : Let Q be the disk bundle of the vector bundle $O_n$(G) constructed in Lemma 8 and let q = dimension Q. By assumption G acts freely on the sphere bundle ∂Q and hence by lemma 3 there is an equivariant embedding of ∂Q into EG(m) x $X_m$ with trivial normal bundle where m = q − 1 − dimension G.

Hence, if a = dimension EG(m) x $X_m$ then there is a tubular neighborhood of the image of ∂Q in EG(m) x $X_m$ equivariantly diffeomorphic to ∂Q × $\mathbb{R}^{a-q+1}$.

Let $Z_n$ = (EG(m) x $X_m$ − ( ∂Q × int($D^{a-q+1}$) ) ∪ (Q × $S^{a-q}$) where we have identified ∂Q with its image in EG(m) x $X_m$. We show that G-manifolds ∈ $\mathcal{O}$(C(G, $\mathcal{O}$)) can be embedded equivariantly in $Z_n$ with a normal bundle subordinate to $\mathbb{R}$.

Let $M^n$ be a G manifold in C(G, $\mathcal{O}$) and let E = ν($M^G$, M) parameterize a tubular neighborhood of $M^G$ in M. Then by lemma 8 we have that the k-dimensional G-vector bundle E which is subordinate to V by definition of a G manifolds in C(G, $\mathcal{O}$) can be equivariantly embedded in $O_n$(G) with a trivial normal bundle. We let i : D(E) → $O_n$(G) be (the restriction of) such an embedding. Since G acts freely on the manifold M - int(D(E)) there is an equivariant embedding of M - int(D(E)) in EG(q) x $X_m$ with trivial normal bundle by Lemma 3. Let j: M - int(D(E)) → EG(m) x $X_m$ be such an embedding. We choose j such that image j is disjoint from image i. Then i|S(E) and j|S(E) are embeddings of S(E) into EG(m) x $X_m$ and hence are isotopic via an isotopy having a trivial normal bundle. We choose the isotopy to avoid the images of i and j. Since S(E) is collared in M we can write M as M - D(E) ∪ S(E) × I ∪ D(E) and produce a map f : M → $Z_n$ using i, the isotopy and j respectively on the 3 submanifolds to get an embedding with a normal bundle that is subordinate to ℝ. Of course we cannot claim that the normal bundle is trivial. However, lemma 11 then yields an embedding of M into $Z_n$ × $X_n$ with a trivial normal bundle.

Next, we must show uniqueness. Let $f_0$ and $f_1$ be equivariant embeddings of M into $Z_n$ × $X_n$ with a trivial normal bundle. We proceed as in the above argument to produce an equivariant embedding of F : M × I → $Z_n$ × $X_n$ with a trivial normal bundle such that f| M × {0} = $f_0$ and f| M × {1} = $f_1$. Thus the map M × I → $Z_n$ × $X_n$ × I , (x,t) → (F(x, t), t) will be the desired isotopy.
In the proof of existence of an embedding we constructed an embedding on D(E) and M – D(E) and then on a collar of S(E). Similarly, we will write M × I as such a union.
Both $f_0$ and $f_1$ map $M^G$ to $(Z_n \times X_n)^G = G_V(V^t)^G \times X_{n-k} \subset D(E)$ and hence there is an invariant tubular neighborhood O of $M^G$ that maps to the interior of D(E) ≈ E. By Lemma 8, $f_0|O$ and $f_1|O$ are isotopic and we denote the isotopy by h: O × I → D(E) × I.
Let k : M → I be a smooth invariant function satisfying $k(M^G) = 0$, k(M – O) = 1 with k increasing radially on O and let K : M × I → I be defined by K(x, t) = 4 t (1-t) k(x). We define F on $K^{-1}$([0, 1/2]) by F(x, t) = ($f_0$(x), t) if t ≤ 1/4, F(x, t) = ($f_1$(x), t) if t ≥ 3/4, F(x, t) = (h(x, 2t–1/2), t) if 1/4 < t < 3/4.

Since G acts freely on $K^{-1}$([3/4, 1]) there is an equivariant embedding L: $K^{-1}$([3/4, 1]) → EG(q) x $X_m$ with trivial normal bundle by Lemma 3 and hence into $Z_n$ × $X_n$ via transversality. Finally, we define F on the collar $K^{-1}$([1/4, 3/4]) using the isotopy between L|$K^{-1}$({3/4}) and F|$K^{-1}$({1/4}).

Thus $Z_n$ × $X_n$ is the desired G-space $U_n$ (𝒸(G, 𝒪)).

Note that the construction is similar to ordinary surgery which uses $\partial(D^a \times D^b) = S^{a-1} \times D^b \cup D^a \times S^{b-1}$; we use $\partial(Q^a \times D^b) = \partial Q \times D^b \cup Q^a \times S^{b-1}$.

Theorem 14 Let If 𝒪 is admissible and |𝒪'| is finite then there exists an n-universal embedding space $U_n$(G, 𝒪) ∈ 𝒸(G, 𝒪) that is, any $M^n$ ∈ 𝒸(G, 𝒪) admits an equivariant embedding in $U_n$(G, 𝒪) with a trivial normal bundle and furthermore any two such embeddings are isotopic.

Proof By induction on on |𝒪'|. If |𝒪'|=1 then 𝒪 consists of an orbit type of the form (H, $\mathbb{R}^m$) for some fixed conjugacy class (H). By lemma 4 we may take $U_n$(G, 𝒪) = G $x_N$ EN'(s + 1) x $X_s$ in that case where N' = N(H)/H and s = dimension M - dimension G/H.

Next assume by induction that theorem 14 holds for all n if |𝒪'| ≤ s -1; we will show it then holds for all n if |𝒪'| = s.

Let $\mathcal{O}$ be an admissible collection of orbit types with $|\mathcal{O}'| = s$. Choose an H among the $[H,W] \in \mathcal{O}'$ such that dimension H is maximal and the number of components of H is maximal among all such subgroups and let $\mathcal{O}^- = \{(K,V) \in \mathcal{O} | K \text{ is not conjugate to H}\}$. Then $\mathcal{O}^-$ is an admissible collection and $|(\mathcal{O}^-)'| = s - 1$ and by our induction hypothesis for any n there is an n-universal embedding space for $\mathfrak{c}(G, \mathcal{O}^-)$, $U_n(G, \mathcal{O}^-)$. We will next construct a $U_n(G, \mathcal{O})$. For each $[H,W] \in \mathcal{O}'$ we have, by Lemma 11, an n-universal embedding space $O_n(G, [H,W])$ which is a G-vector bundle so the disjoint union of all $O_n(G, [H,W])$ with $[H,W] \in \mathcal{O}'$ is a G-vector bundle (of possibly varying fiber dimension). The sphere bundle of this disjoint union is a G-manifold in $\mathfrak{c}(G, \mathcal{O}^-)$ and hence admits an equivariant embedding into $U_m(G, \mathcal{O}^-)$ for some m.

We use this embedding to attach the disk bundle to $U_m(G, \mathcal{O}^-)$ as in Propsition 13 to get a G-manifold $Z_n$ such that every manifold in $\mathfrak{c}(G, \mathcal{O})$ can be embedded in $Z_n$ with bormal bundle subordinate to $\mathbb{R}$. Finally, we invoke lemma 11 to see that $Z_n \times X_n$ is the desired G-space $U_n(\mathfrak{c}(G, \mathcal{O}))$.

Let $\mathfrak{c}^\mathfrak{c}(G, \mathcal{O})$ denote $\{M^n \in \mathfrak{c}(G, \mathcal{O}) | M \text{ is compact}\}$.

Corollary 15 If $\mathcal{O}$ is admissible and $|\mathcal{O}'|$ is finite then there exists an n-universal embedding space $U_n^\mathfrak{c}(G, \mathcal{O}) \in \mathfrak{c}^\mathfrak{c}(G, \mathcal{O})$ such that any $M^n \in \mathfrak{c}^\mathfrak{c}(G, \mathcal{O})$ admits an equivariant embedding in $U_n^\mathfrak{c}(G, \mathcal{O})$ with a trivial normal bundle and furthermore any two such embeddings are isotopic.

Proof: The construction in theoem 14 is compact if $|\mathcal{O}'|$ is finite.

For G finite, theorem 14 is sufficient. However, the condition that $|\mathcal{O}|$ is finite will not be satisfied in general for an arbitrary compact Lie group G.

There are several obstacles to eliminating this condition. First of all the simple induction argument does not work if $|\mathcal{O}'|$ is not finite. Here is how we remedy that.

Definitions If W is a representation of a group H we denote by $\widehat{W}$ the representation $W/W^H$.
We define an order relation $\prec$ on $\mathcal{O}'$ by $[H,W] \prec [H', W']$ if there is a map j in $\mathfrak{c}(G, \mathcal{O})$,

$j: G \times_H \widehat{W} \to G \times_{H'} \widehat{W'}$.

We note that $\prec$ is reflexive and transitive and that $[H,W] \prec [H, W \oplus \mathbb{R}]$ and $[H, W \oplus \mathbb{R}] \prec [H, W]$.

Lemma 16 If $[H,W] \prec [H', W']$ and $[H,W] \neq [H', W']$ then dim $W <$ dim $W'$.

Proof. Since j is an embedding dimension $G \times_H \widehat{W} \leq G \times_{H'} \widehat{W'}$ so dim $G -$ dim $H +$ dim $W \leq$ dim $G -$ dim $H' +$ dim $W'$ or dim $W' -$ dim $W \geq$ dim $H -$ dim $H'$ and since $H = G_{[j(g, 0])}$is conjugate to a subgroup of H' dim $H \leq$ dim $H'$ thus dim $W' -$ dim $W \geq 0$. We assume dim $W' -$ dim $W = 0$ and show a contradiction. If $j([e, 0]) = [e, 0]$ then $H = H'$ and $dj_{[e,0]}$ maps the slice at $[e, 0]$, $W$, monomorphically to the slice at $j([e.0]) = W'$ since $j$ is an embedding. Since the normal bundle of $j$ is trivial by assumption and $W'^H = 0$, the cokernel of $dj_{[e,0]} = 0$ and thus $dj$ is an isomorphism and $[H, W] = [H', W']$ which is a contradiction. So assume $j([e, 0]) = [e, x]$ with $x \neq 0$; then $H$ is a proper subgroup of $H'$ and $dj_{[e,0]}$ maps the slice at $[e, 0]$, $W$, monomorphically to the slice at $j([e.0]) = [e,x]$ which is $W'|H$ since $\dim H = \dim H'$. Since $0 \neq x \in W'^H$, $dj(W) \neq W'$ so $\dim W < \dim W'$. Finally, if $j([e, 0]) = [g, x]$ we note that the ector space $W'' = g \times W'$ is a representation of the group $H'' = gH' g^{-1}$ and $[H', W'] = [H'', W'']$ so the previous argument applies and $\dim W < \dim W'' = \dim W'$.

A consequence of Lemma 16 is that for any chain of inequalities $[H_0, 0] \prec [H_1, W_1]...[H_s, W_s]$,

$s \leq \dim \widetilde{W}$.

We next define a function $l$: $\mathcal{O}' \to \mathbb{N}$ **by** $l([H, 0]) = 0$ **and** $l([H, W]) = s$ where $s$ is the maximal length of a chain of inequalities $[H_0, 0] \prec [H_1, W_1] \prec .. \prec [H_{s-1}, W_{s-1}] \prec [H, W]$. As noted above $l$ is well defined.

Lemma 16 shows that $l$: $\mathcal{O}' \to \mathbb{N}$ is well defined even **if** $|\mathcal{O}'|$ is not finite. Furthermore, if $M$ is an n-dimensional G-manifold then for every $x \in M$, $l([G_x, S_x]) \leq n$ for every pair $G_x, S_x$ where $G_x$ is the isotropy subgroup at x and $S_x$ is the slice representation of $G_x$ at x

With proposition 14 illustrating the method of construction of n-universal embedding spaces we are ready to prove theorem 17 by induction using $l$.

Theorem 17 If $\mathcal{O}$ is an admissible collection of orbit types then there exists an n-universal embedding space $U_n(G, \mathcal{O}) \in \mathcal{C}(G, \mathcal{O})$ that is, any $M^n \in \mathcal{C}(G, \mathcal{O})$ admits an equivariant embedding in $U_n(G, \mathcal{O})$ with a trivial normal bundle and furthermore any two such embeddings are isotopic.

Proof Let $\mathcal{O}[a] = \{ [H,W] \in \mathcal{O}' | l([H,W]) \leq a \}$ and let $\mathcal{C}(G, \mathcal{O}[a]) = \{ M \in \mathcal{C}(G, \mathcal{O}) | $ for all $x \in M$ $[G_x, S_x] \in \mathcal{O}[a]\}$.

Fix n. We prove the existence of $U_n(G, \mathcal{O}[a])$ by induction on a. We start with $a = 0$ and prove the existence of $U_t(G, \mathcal{O}[0])$ for all t.

Let $B = \amalg\, G \times_{N(H)} EN'(H)(n + 1) \times X_t$ be the disjoint union over all $[H, 0]$ in $\mathcal{O}'$) such that $l([H, \mathbb{R}]) = 0$. It is a union of manifolds of possibly varying dimension: dimension $G \times_{N(H)} EN'(H)(n + 1) \times X_n =$ dimension $G/N(H) +$ dimension $X_n +$ dimension $EN'(H)(n + 1)$ By lemma 4 any $M^n \in \mathcal{C}(G, \mathcal{O}[0])$ embeds in B with trivial normal bundle and any two such embeddings are equivariantly isotopic.

Remark: The distinct spaces $G \times_{N(H)} EN'(H)(n + 1) \times X_n$ above correspond to distinct principal isotropy subgroups and $U_n(G, \mathcal{O})/G$ will necessarily have have as many components as $G$ has conjugacy classes of subgroups.

There is a problem: the dimensions of the components of B may not be bounded so the padding trick will not work. The remedy for this problems is simple:

Lemma 18 Modified Steenrod's construction:

Let G be a compact Lie group of dimension d and let $n \geq 0$; then there is a smooth manifold E'G(n) with a free action of G such that if G acts freely on a manifold $M^{n+d}$, B is an invariant submanifold and $h : B \to E'G(n)$ is a smooth equivariant map then there is a smooth equivariant map f from M to E'G(n) such that $f|B = h$. Furthermore dimension $E'G(n) \leq 2n + d + 1$.

Proof Let $\pi : EG(n) \to EG(n)/G$ be the principal bundle constructed by Steenrod. To construct E'G(n) we note that EG(n)/G is triangulable; let K denote the n-dimensional skeleton of EG(n)/G in some triangulation. We first note that $\pi : \pi^{-1}(K) \to K$ is also n-universal. If G acts freely on an (n+d)-dimensional manifold $M^{n+d}$, B is an invariant submanifold and $h : B \to \pi^{-1}(K) \subset EG(n)$ is a continuous equivariant map then h extends to a continuous equivariant map $H: M \to EG(n)$ and the quotient map $H' : M/G \to EG(n)/G$ is homotopic to a map $j: M/G \to K$ by the simplicial approximation theorem. Thus by the covering homotopy theorem H is homotopic to a map $J: M \to \pi^{-1}(K)$.

However, the n-universal G-space $\pi^{-1}(K)$ is not a manifold. To remedy that we embed K in $\mathbb{R}^{2n+1}$ and let $\Omega$ be a regular neighborhood of K, that is, $\Omega$ is a smooth manifold with boundary and K is a deformation retract of $\Omega$. [RS] Let $r : \Omega \to K$ be the retraction and $r" : \Omega \to EG(n)/G$ be a smooth map homotopic to $r' : \Omega \to K \to EG(n)/G$. We set $E'G(n) = r"^*EG(n)|K$ and claim that this G-manifold has the desired properties.

We note first that the free G-manifold E'G(n) has dimension $E'G(n) = d + 2n + 1$.

Also, since K is a deformation retract of $\Omega$ the principal bundles $\pi : \pi^{-1}(K) \to K$ and $E'G(n) \to \Omega$ both have the n-universal property. Thus E'G(n) is the desired modification of the Steenrod construction.

We now resume the proof of Theorem 17.
Let $B' = \amalg\ G \times_{N(H)} E'N'(H)(n + 1) \times X_n$ be the disjoint union over all $[H, 0]$ in $\mathcal{O}'$) such that $l([H, \mathbb{R}]) = 0$. It is a union of manifolds of possibly varying dimension:
dimension $G \times_{N(H)} EN'(H)(n + 1) \times X_n =$ dimension $G/N(H) +$ dimension $X_n +$ dimension $EN'(H)(n + 1)$
$\leq d +$ dimension $X_n +$ dimension $N'(H) + 2n + 3 \leq 2d +$ dimension $X_n + 2n + 3$ which is independent of H so we can pad where necessary so that each of the summands has the same dimension. That modification of B' yields $U_t(G, \mathcal{O}[0])$, the desired G-manifold.

Now suppose inductively that we have constructed a $U_t(G, \mathcal{O}[a])$ for all t. We must show the existence of a $U_n(G, \mathcal{O}[a+1])$.

In the proof of Theorem 14 we did induction by surgically attaching one disk bundle $O_n(G, [H,W])$ at a time; that was possible since $|\mathcal{O}|$ was finite. In the present case we will surgically attach all disk bundles $O_n(G, [H,W])$ where $l(H) = a+1$. (The bundles $O_n(G, [H,W])$ are defined and constructed in Lemma 11. The bundles $O_n(G, [H,W])$ are of the unions of G-vector bundles of the form $G \times_N (\mu_k | G_k(tV)^H \times EN'(s+1)) \times X_s$ where $s = n - k - $ dimension $G/H$ and V is a representation of G that depends on H.)

Since the number of pairs [H, W] may be infinte the dimensions $O_n(G, [H,W])$ may not be bounded even if we replace $EN'(s+1))$ above by $E'N'(s+1))$; the problem is that the dimension of $G_k(tV)^H$ depends on dimension V which is not obviously bounded. Hwever, the remedy to this problem is similar.

Each $G \times_N (\mu_k | G_k(tV)^H \times EN'(s+1)) \times X_s$ is a G-vector bundle over $G \times_N G_k(tV)^H \times EN'(s+1)) \times X_s$ which is a G/H bundle over $(G_k(tV)^H \times EN'(s+1)) \times X_s)/N$ which is a smooth manifold. Let K denote the s + 1 dimensional skeleton of $(G_k(tV)^H \times EN'(s+1)) \times X_s)/N$ in some triangulation where $s = n - k - $ dimension $G/H$. K is not a manifold so we embed K in $\mathbb{R}^{2s+1}$ and let $\Omega$ be a regular neighborhood of K, that is, $\Omega$ is a smooth manifold with boundary and K is a deformation retract of $\Omega$. [RS] Let $r : \Omega \to K$ be the retraction and $r'' :$
$\Omega \to (G_k(tV)^H \times EN'(s+1)) \times X_s)/N$ be a smooth map homotopic to
$r' : \Omega \to K \to (G_k(tV)^H \times EN'(s+1)) \times X_s)/N$.
We set $O_n(G, [H,W])' = r''^*(G \times_N (\mu_k | G_k(tV)^H \times EN'(s+1)) \times X_s)$ and claim that, as in the proof of $U_t(G, \mathcal{O}[0])$, existence of this G-vector bundle of bounded dimension has the desired properties. We now wish to surgically attach the disk bundles of all $O_n(G, [H,W])'$ where $l(H) = a+1$ as in Proposition 13 and Theorem 14. The bundles $D(Q)$ attached in the proof of Proposition 13 and $D(O_n(G, [H,W]))$ attached in the the proof of Theorem 14 had base spaces that were closed manifolds. The bundles $O_n(G, [H,W])'$ howver, do not have closed manifolds as base space. To avoid discussion of manifolds with corners we simply extend our bundles to the double of $\Omega$ which is a closed manifold, attach as before and then remove that part of the disk bundle over one copy of $\Omega$. That completes the construction of $U_n(G, \mathcal{O}[a+1])$ and the inductive step in the proof.

We note that there is no possibility of constructing a compact version of $U_n(G, \mathcal{O})$ as in Corollary 15 when $|\mathcal{O}|$ is not finite. In general $.U_n(G, \mathcal{O})$ will have an infinite number of components.

We now proceed to some applications of Theorem 17.

§3 Applications

In this section $\mathcal{O}$ will denote an admissible collection of orbit types. Let $\mathfrak{c}(G, \mathcal{O})$ be the category of G

manifolds with $(G_x, S_x) \in \mathcal{O}$ for all $x \in \mathcal{O}$ and with morphisms being equivariant embeddings with trivial normal bundles and and for each n let $U_n(\mathsf{C}(G, \mathcal{O}))$ be an n-universal embedding space, that is, for every n-dimensional manifold M in $\mathsf{C}(G, \mathcal{O})$ there is an equivariant embedding of M into $U_n(\mathsf{C}(G, \mathcal{O}))$ with a trivial normal bundle. A specific trivialization, i.e., an equivariant bundle map
$T : \nu(M, U_n(\mathsf{C}(G, \mathcal{O}))) \to M \times \mathbb{R}^N$ where dimension $U_n(\mathsf{C}(G, \mathcal{O}))$ is $n + N$ is called a framing. In general, distinct framings need not be homotopic–consider the trivial example of embeddings of zero dimensional manifolds into $\mathbb{R}^2$. However we have

**Proposition 19** If $M \in \mathsf{C}(G, \mathcal{O})$ is an n-dimensional manifold and $T_1$ and $T_2$ are framings of $\nu(M, U_{2n+1}(\mathsf{C}(G, \mathcal{O})))$ for some choice of $U_{2n+1}(\mathsf{C}(G, \mathcal{O}))$ then $T_1$ and $T_2$ are equivariantly homotopic.

**Proof** Let $h : M \to GL(N, \mathbb{R})$ be the map defined by $T_1 \circ T_2^{-1} : M \times \mathbb{R}^N \to M \times \mathbb{R}^N$
The map h is homotopic **to** a map $h' : M \to GL(n, \mathbb{R})$.
(The map $GL(n, \mathbb{R}) \to GL(N, \mathbb{R})$ is $n - connected$.).
Let $L = \left(M \times \mathbb{R}^N \times \mathbb{R}\right) / (x, v, t) \sim (x, h'(v), t + 1)$.

Since dimension $L = 2n + 1$ there is an embedding of L in $U_{2n+1}(C(G, \mathcal{O}))$ with trivial normal bundle. The embedding of M in L, $x \to (x, 0, 0)$ is isotopic to the embedding of M in L, $x \to (x, 0, 1)$. Thus the trivializations of M in $U_{2n+1}(C(G, \mathcal{O}))$ are homotopic.

We now define the final version of the definition of an n-universal embedding space:

**Definition** A space $Y_n(\mathsf{C}(G, \mathcal{O})) \in \mathsf{C}(G, \mathcal{O})$ is an n-universal embedding space for the category $\mathsf{C}(G, \mathcal{O})$ if any n-dimenional manifold $M \in \mathsf{C}(G, \mathcal{O})$ can be equivariantly embedded in $Y_n(\mathsf{C}(G, \mathcal{O}))$ with a trivial normal bundle and any 2 such embeddings are unique up to equivariant isotopy that is, if j0 and j1 are embeddings of $M \subset Y_n(\mathsf{C}(G, \mathcal{O}))$ and J0: $\nu(M, Y_n(\mathsf{C}(G, \mathcal{O}))) \to M \times \mathbb{R}^N$,
J1: $\nu(M, Y_n(\mathsf{C}(G, \mathcal{O}))) \to M \times \mathbb{R}^N$ are trivializations of the normal bundles of the respective embeddings than there is an embedding j of $M \times I \subset Y_n(\mathsf{C}(G, \mathcal{O})) \times I$ and a trivialization of

$J : \nu(M \times I, Y_n(\mathsf{C}(G, \mathcal{O})) \times I) \to M \times I \times \mathbb{R}^N$ such that $j | M \times \{0\} = j0$,
$j | M \times \{1\} = j1$, $J | \nu(M \times I, Y_n(\mathsf{C}(G, \mathcal{O})) \times I) | M \times \{0\} \times \mathbb{R}^N = J0$, **and**
$J | \nu(M \times I, Y_n(\mathsf{C}(G, \mathcal{O})) \times I) | M \times \{1\} \times \mathbb{R}^N = J1$.

In view of the above proposition the constructions of section 2 all yield $n -$ universal embedding spaces.

**Corollary 20** If $Y_{n+1}(\mathsf{C}(G, \mathcal{O}))$ is an (n+1)-universal embedding space for the category $\mathsf{C}(G, \mathcal{O})$ and dimenson $Y_{n+1}(\mathsf{C}(G, \mathcal{O})) = n + N$ then $[Y_{n+1}(\mathsf{C}(G, \mathcal{O}))^*, \infty ; (X \times D^N)/(X \times S^{N-1})]_G =$

$\mathscr{U}_n(X, C(G, \mathcal{O}))$ where * indicates the 1-point compactification and $[A, B]_G$ denotes the G-equivariant homotopy classes of maps $A \to B$ and $\mathscr{U}_n(X, C(G, \mathcal{O}))$ denotes the G−equivariant bordism group of the space X.

In particular, if X = pt then $\mathscr{U}_n(pt, C(G, \mathcal{O})) = [Y_{n+1}(C(G))^*/G, \infty; S^N]$, i.e. the equivariant bordism groups are the ordinary cohomotopy groups of the space $Y_{n+1}(C(G))^*/G$.

Proof. Let $\Phi : [Y_{n+1}(C(G))^*, \infty; (X \times D^N)/(X \times S^{N-1})]_G \to \mathscr{U}_n(X, C(G))$ be given by $\Phi(F) = [F'^{-1}(X \times \{0\}), F'|F'^{-1}(X \times \{0\})]$ where $F'$ is equivariantly homotopic to $F$ and transverse regular to $X \times \{0\}$. This yields a G-manifold and an equivariant map to X. The usual argument shows that homotopic maps yield bordant manifolds and maps so $\Phi$ is well defined. Given [M, f] in $\mathscr{U}_n(X, C(G, \mathcal{O}))$ one embedds M in $Y_{n+1}(C(G, \mathcal{O}))$ with trivial normal bundle and performs the Thom-Pontryagin construction to see that $\Phi$ is onto. [T] To see that $\Phi$ is one to one suppose that $\Phi(F_0) = [M_0, f_0]$ and $\Phi(F_1) = [M_1, f_1]$ and that there is a compact manifold with boundary W in $C(G, \mathcal{O})$ and a map $f : W \to X$ such that $\partial W = M_0 \cup M_1$ and $f_0 = f|M_0, f_1 = f|M_1$. Let h: $W \to [0, 1]$ be a smooth equivariant map with $h^{-1}(0) = M_0, h^{-1}(1) = M_1$. Since dimension double(W) is n+1 there is an equivariant embedding k: double(W) $\to Y_{n+1}(C(G, \mathcal{O}))$ with trivial normal bundle and hence an embedding (k, h) : W $\to Y_{n+1}(C(G, \mathcal{O})) \times [0, 1]$ with trivial normal bundle. We then perform the Thom-Pontryagin construction on W $\subset Y_{n+1}(C(G, \mathcal{O})) \times [0, 1]$, $f : W \to X$ to get an equivariant homotopy $H : Y_{n+1}(C(G, \mathcal{O}))^* \times [0,1], \infty; (X \times D^N)/(X \times S^{N-1})$ with $(H| Y_{n+1}(C(G, \mathcal{O}))^* \times \{0\})^{-1} \{(X \times S^{N-1})\} = M_0$ and $(H| Y_{n+1}(C(G, \mathcal{O}))^* \times \{1\})^{-1} \{(X \times S^{N-1})\} = M_1$ and $f_0 = f|M_0, f_1 = f|M_1$. The trivializations of
$\nu(M_0, Y_{n+1}(C(G, \mathcal{O}))$ given by $F_0$ and $H| Y_{n+1}(C(G, \mathcal{O})) \times \{0\})$ are homotopic by the above Proposition hence $H| Y_{n+1}(C(G, \mathcal{O}))^* \times \{0\}$ is homotopic to $F_0$ and similarly, $H| Y_{n+1}(C(G, \mathcal{O}))^* \times \{1\}$ is homotopic to $F_1$ thus $F_0$ is homotopic to $F_1$.

We have shown that $\Phi$ is an isomorphism of sets. However, $[Y_{n+1}(C(G, \mathcal{O}))^*, \infty; (X \times D^N)/(X \times S^{N-1})]_G$
has a group structure if n+1 < N since $(X \times D^N)/(X \times S^{N-1}) = S^N \wedge X_+$ [A2] and
$S^N \wedge X_+ \vee S^N \wedge X_+ \to S^N \wedge X_+ \times S^N \wedge X_+$ is an equivariant cofibration **and**
$S^N \wedge X_+ \wedge S^N \wedge X_+ = S^{2N} \wedge X_+$ is N connected and thus, any map of $Y_{n+1}(C(G, \mathcal{O}))$ can be homotoped into $S^N \wedge X_+ \vee S^N \wedge X_+$ and the folding map as usual gives a group structure to the equivariant cohomotopy group. We conclude that $\Phi$ is a group isomorphism.

Note that in the case that X = pt (or more generally if $X = X^G$) the group action on $S^N \wedge X$ is trivial so any equivariant map from $Y_{n+1}(C(G, \mathcal{O}))$ to $S^N \wedge X$ must factor throught $Y_{n+1}(C(G, \mathcal{O}))/G$.

Remark The cases $X = X^G$ and X a free G-space should be amenable to further simplification.

Let h*( ) be a cohomology theory defined for all manifolds M in $C(G, \mathcal{O})$. An h* natural class for $C(G, \mathcal{O})$ is a cohomology class $z(M) \in h^*(M)$ defined for all manifolds M in $C(G, \mathcal{O})$ such that if $f : A \to B$ is a map in $C(G, \mathcal{O})$ then $z(A) = f^*(z(B))$.

We denote by $Nat(C(G, \mathcal{O}), h^*)$ the set of natural classes for h*.

Note The only property of h* that we use is that it is a contravariant functor.

When G is the trivial group natural classes are just the h* groups of the classifying space for the tangent bundles of our manifolds. Things are much more interesting for nontrivial G. For example, fixed point sets of subgroups give rise to natural classes via Poincaré duality. [W2]

The computation of $Nat(C(G, \mathcal{O}), h^*)$ might use the following:

Theorem 21 $Nat(C(G, \mathcal{O}), h^*) = \text{inv}\lim_{n \to \infty} (h^*(U_n(C(G, \mathcal{O}))))$.

Proof. Let z be a natural class for $C(G, \mathcal{O})$; then for any n-universal embedding spaces $U_n(C(G, \mathcal{O}))$ for $C(G, \mathcal{O})$, $z(U_n(C(G, \mathcal{O})))$ is defined since $U_n(C(G, \mathcal{O})) \in C(G, \mathcal{O})$. Thus we get an element $\in$ $\text{inv}\lim_{n \to \infty} (h^*(U_n(C(G, \mathcal{O}))))$.

Now let $z \in \text{inv}\lim_{n \to \infty} (h^*(U_n(C(G, \mathcal{O}))))$, that is, a set of elements $x_n, \in h^*(U_n(C(G, \mathcal{O})))$ and let $M \in C(G, \mathcal{O})$. There is an equivariant embedding $i : M \subset U_n(C(G, \mathcal{O}))$ for some n by theorem 17. Define $z(M) = i^* x_n$. Note z(M) does not depend on which specific i is chosen since the choice is unique up to equivariant isotopy. Furthermore, if $j : M \subset U_m(C(G, \mathcal{O}))$ is another embedding of M in an m-universal embedding space, we choose $s > \dim U_n(C(G, \mathcal{O}), \dim U_m(C(G, \mathcal{O}))$ and embeddings $f1 : U_n(C(G, \mathcal{O}) \subset U_s(C(G, \mathcal{O}))$, $f2 : U_m(C(G, \mathcal{O}) \subset U_s(C(G, \mathcal{O}))$ and then note that $f1 \circ i$ and $f2 \circ j$ and both embeddings of $M \subset U_s(C(G, \mathcal{O}))$ and hence are equivariantly isotopic. In particular, $(f1 \circ i)^*(x_s) = (f2 \circ j)^*(x_s)$ so $i^* \circ f1^*(x_s) = i^*(x_n) = j^* \circ f2^*(x_s) = j^*(x_m)$; thus z(M) does not depend on the choice of embedding space. Finally, if $f : A \to B$ is an embedding and $i : B \subset U_n(C(G, \mathcal{O}))$ then $z(B) = i^*(x_n)$, $f^*(z(B)) = f^*(i^*(x_n)) = (i \circ f)^*(x_n) = z(A)$ where the last equality follows from the fact that $i \circ f$ is also an embedding of A in an n-universal embedding space.

Suggestions for further investigation:
a) calulate $\mathcal{N}_n(pt, C(G, \mathcal{O}))$ for some interesting G
b) construct universal embedding spaces for other types of G-manifolds, for example, oriented or stably almost complex manifolds

c) calculate $Nat(C(G, \mathcal{O}), h^*)$ insome nontrivial cases. A rudimentary start to this problem is {W2}.

Department of Mathematics University of Michigan Ann Arbor, Michigan 48109

*Email address:* awass@umich.edu